\documentclass{birkjour_modified}

\usepackage{amsfonts}
\usepackage{amsmath}
\usepackage{bbold}
\usepackage{amsthm}
\usepackage{bbm}
\usepackage{graphicx}
\usepackage{subcaption}
\usepackage[margin=0.5cm]{caption}

\def\E{\mathcal E}
\def\B{\mathcal B}
\renewcommand{\vec}[1]{\mathbf{#1}}
\DeclareMathOperator{\tr}{tr}

\newtheorem{thm}{Theorem}
\newtheorem{cor}[thm]{Corollary}

\theoremstyle{definition}

\theoremstyle{remark}

\begin{document}

\title{The Euler and Grace-Danielsson inequalities for nested triangles and tetrahedra:\\ a derivation and generalisation using \\quantum information theory}

\author{Antony Milne}
\address{Controlled Quantum Dynamics Theory, Department of Physics,\\ Imperial College London, London SW7 2AZ, UK \\ E-mail: \texttt{antony.milne@gmail.com}\\}

\begin{abstract}We derive several results in classical Euclidean elementary geometry using the steering ellipsoid formalism from quantum mechanics. This gives a physically motivated derivation of very non-trivial geometric results, some of which are entirely new. We consider a sphere of radius $r$ contained inside another sphere of radius $R$, with the sphere centres separated by distance $d$. When does there exist a `nested' tetrahedron circumscribed about the smaller sphere and inscribed in the larger? We derive the Grace-Danielsson inequality $d^2\leq (R+r)(R-3r)$ as the sole necessary and sufficient condition for the existence of a nested tetrahedron. Our method also gives the condition $d^2\leq R(R-2r)$ for the existence of a nested triangle in the analogous 2-dimensional scenario. These results imply the Euler inequality in 2 and 3 dimensions. Furthermore, we formulate a new inequality that applies to the more general case of ellipses and ellipsoids.\end{abstract}

\maketitle

\section{Introduction}

Let $d$ be the distance between the circumcentre and incentre of a triangle with circumradius $R$ and inradius $r$. Independently, Chapple (in 1746) and Euler (in 1765) found a classic result of 2-dimensional Euclidean geometry that holds for all triangles: $d^2=R(R-2r)$~\cite{GeomIneq}. The analogous scenario in 3-dimensional space involves a tetrahedron, circumsphere and insphere. In 1816, Gergonne asked whether in 3 dimensions $d$ could be similarly expressed as a function of only $R$ and $r$. Eight years later, Durrande gave the solution $d^2=(R+r)(R-3r)$. This was widely accepted for many years but is in fact incorrect, and there cannot exist such an equality that holds for all tetrahedra (see Ref.~\cite{Enigma} for a full discussion).

\begin{figure}
\centering
\begin{subfigure}{.5\textwidth}
  \flushleft
  \includegraphics[width=0.95\linewidth]{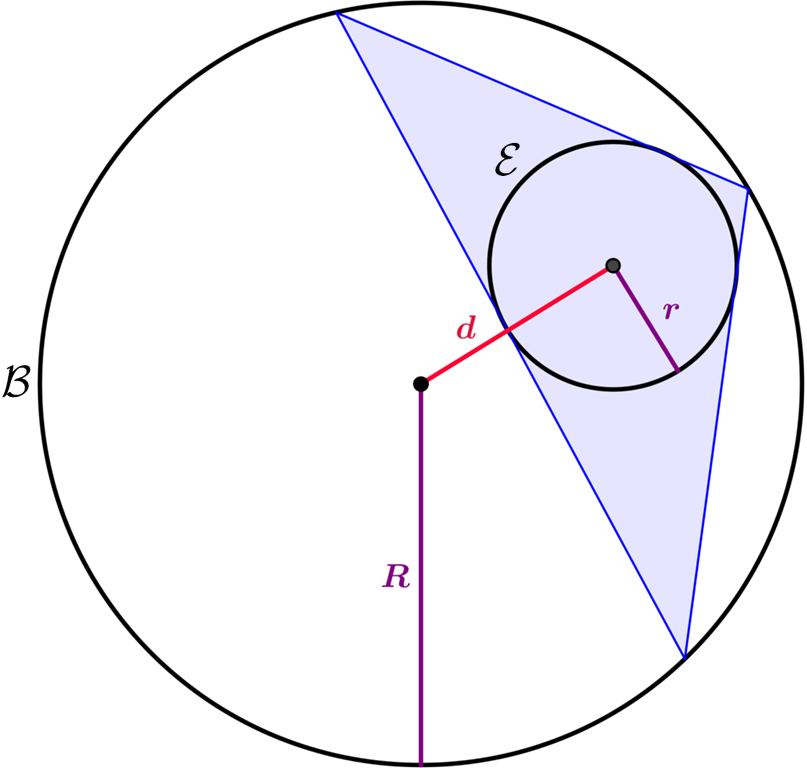}
  \caption{A nested triangle}
\end{subfigure}%
\begin{subfigure}{.5\textwidth}
	\flushright
    \includegraphics[width=0.95\linewidth]{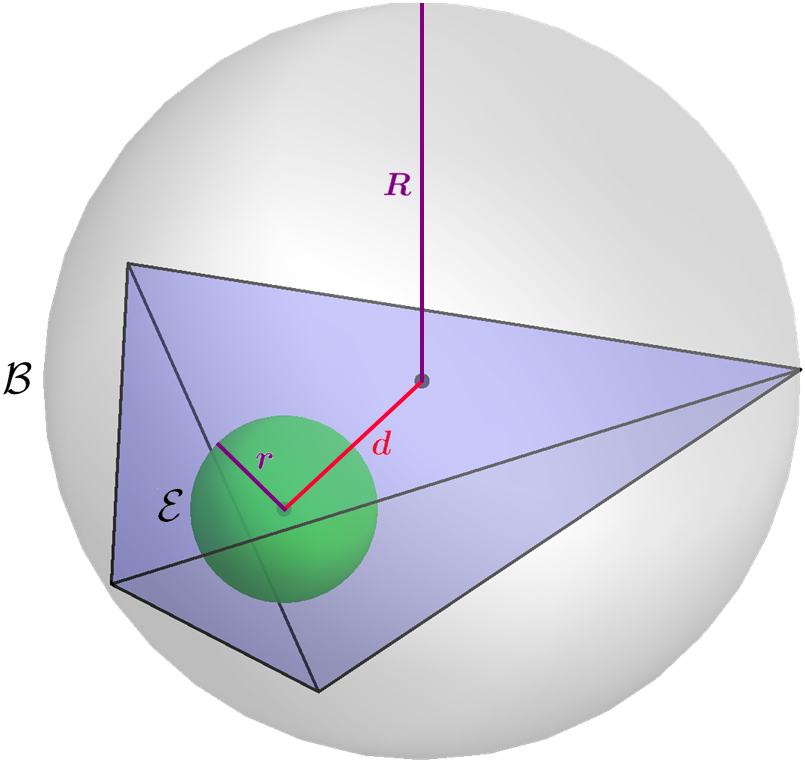}
  \caption{A nested tetrahedron}
\end{subfigure}
\caption{A circle (sphere) $\E$ of radius $r$ is contained inside another circle (sphere) $\B$ of radius $R$. The distance between the centres of $\E$ and $\B$ is $d$. In the examples shown here there exists a nested triangle (tetrahedron) circumscribed about $\E$ and inscribed in $\B$.}
\label{fig:nested}
\end{figure}

We consider a closely related question. A circle (sphere) $\E$ of radius $r$ is contained inside another circle (sphere) $\B$ of radius $R$. Let the distance between the centres of $\E$ and $\B$ be $d$. What are the necessary and sufficient conditions for the existence of a triangle (tetrahedron) circumscribed about $\E$ and inscribed in $\B$? We shall describe such a triangle (tetrahedron) as \emph{nested}. Examples of the 2- and 3-dimensional scenarios are shown in Figure \ref{fig:nested}.

In 2 dimensions, the sole condition for the existence of a nested triangle is~\cite{RecentAdvances}\begin{equation}\label{2d_ineq}d^2\leq R(R-2r).\end{equation} In 3 dimensions, Grace (in 1917, see Ref.~\cite{Grace}) and Danielsson (in 1949, see Ref.~\cite{Danielsson}) proved that the sole condition for the existence of a nested tetrahedron is \begin{equation}\label{3d_ineq}d^2\leq (R+r)(R-3r).\end{equation} We shall call \eqref{3d_ineq} the \emph{Grace-Danielsson inequality}.

Ref.~\cite{Enigma} notes that the statement of this problem is formulated in terms of classical Euclidean elementary geometry (CEEG) but Danielsson's proof is based on some intricate projective geometry. This poses a challenge to prove inequality \eqref{3d_ineq} using only methods belonging to CEEG. In this article we prove the Grace-Danielsson inequality without the explicit use of any geometry at all; instead we use a tool from quantum mechanics called the \textit{steering ellipsoid}~\cite{QSE}.

We begin by outlining the theory behind steering ellipsoids. By understanding the separability of two-qubit states in the formalism, we derive a single key inequality that applies to both the 2- and 3-dimensional scenarios (Theorem \ref{keyinequality}). From this inequality the results \eqref{2d_ineq} and \eqref{3d_ineq} easily follow (Corollaries \ref{circle} and \ref{sphere} respectively). Moreover, our inequality can be used to give a necessary and sufficient condition for the existence of a nested triangle (tetrahedron) in the general case that $\E$ is an ellipse (ellipsoid). This is the first time that such a condition has been formulated. We give an example of how it might be used for ellipsoids in Corollary \ref{ellipsoid}.

The $n$-dimensional Euler inequality gives $R \geq nr$ for the circumradius and inradius of an $n$-dimensional simplex~\cite{RecentAdvances}. Since $d^2\geq 0$, inequalities \eqref{2d_ineq} and \eqref{3d_ineq} imply the Euler inequality for $n=2$ and $n=3$ respectively. We will therefore recover these classic results of 2- and 3-dimensional Euclidean geometry through a physical argument. There already exist some physically motivated derivations of geometric results; for example, the generalised parallel axes theorem can be used to prove properties of triangle centres such as Euler's formula $d^2=R(R-2r)$~\cite{GPAT}. Our derivation is particularly remarkable since it gives a new generalisation as well as recovering a result that remains unproven using CEEG.

Finally, we note that this work has led to renewed interest in the $n$-dimensional scenario, with Egan conjecturing $d^2\leq (R+(n-2)r)(R-nr)$ to be necessary and sufficient for the existence of a nested simplex~\cite{Egan}. Although sufficiency of this condition has been proven, there is not yet a proof of its necessity.

\vspace{-3mm}
\section{Quantum steering ellipsoids}

We begin by briefly describing the steering ellipsoid formalism and its significance within the broad field of quantum information theory.

Quantum information theory studies how quantum mechanics can be exploited to process information in ways that cannot be achieved using classical physics~\cite{NC}. The basic unit is a \textit{qubit}, the quantum mechanical analogue of a classical bit of data. A central concept in quantum information is the phenomenon of entanglement. This describes a form of non-classical correlation between physically separated systems and has for a long time been a major curiosity in the foundations of physics. More recently, entanglement has been recognised as a powerful resource for many quantum information processing tasks, including quantum computing.

Any quantum system is associated with a complex Hilbert space $\mathcal{H}$ (a complete vector space equipped with an inner product). The \textit{state} of the quantum system is described by a linear operator $\rho$ acting on $\mathcal{H}$. By definition, $\rho$ must be unit trace and positive semi-definite ($\rho \geq 0$). A single qubit corresponds to $\mathcal{H}=\mathbbm{C}^2$; the associated state $\rho$ is given by a $2\times 2$ Hermitian matrix. We can express this as $\rho=\frac{1}{2}(\mathbb{1}+\vec{r}\cdot \boldsymbol\sigma)$ in the Pauli basis $(\mathbb{1}, \boldsymbol\sigma)$, where \begin{equation}\nonumber\mathbb{1}=\begin{pmatrix} 1&0 \\ 0&1\end{pmatrix},\,\boldsymbol\sigma=\left(\begin{pmatrix} 0&1 \\ 1&0\end{pmatrix},\,\begin{pmatrix} 0&-i \\ i&0\end{pmatrix},\,\begin{pmatrix} 1&0 \\ 0&-1\end{pmatrix}\right)\end{equation}
and $\vec r$ is a real 3-dimensional \textit{Bloch vector}. The constraint $\rho\geq 0$ is equivalent to the condition $|\vec r|\leq 1$. Thus any single qubit state can be represented by a Bloch vector inside the unit sphere (\textit{Bloch sphere}), which we denote $\widetilde\B$.

Two qubits are fundamentally important for studying quantum correlations as they form the simplest system that can exhibit entanglement. A two-qubit system has a state $\rho$ acting on $\mathcal{H}=\mathbbm{C}^2\otimes\mathbbm{C}^2$, where $\otimes$ is the tensor product. In this case $\rho$ corresponds to a $4 \times 4$ Hermitian matrix that can be described using 15 real parameters. The steering ellipsoid formalism was developed as a method for representing any such two-qubit state using a set of 3-dimensional objects that can be easily visualised~\cite{QSE}.

For our purposes it suffices to focus on a restricted set of so-called \textit{canonical} two-qubit states. If we expand in the product Pauli basis $(\mathbb{1}, \boldsymbol\sigma)\otimes(\mathbb{1}, \boldsymbol\sigma)$, these are given by\begin{equation}\label{canonical}\rho=\frac{1}{4}(\mathbb{1}\otimes\mathbb{1}+\vec{d}\cdot \boldsymbol\sigma\otimes\mathbb{1}+\sum_{i,j=1}^3 S_{ij}\sigma_i\otimes\sigma_j).\end{equation}
Say that Alice and Bob each hold a single qubit and are jointly described by this state. The Bloch vector representing Alice's qubit is $\vec d$ and the Bloch vector representing Bob's qubit is $\vec 0$. However, the two Bloch vectors alone do not fully describe $\rho$: there is also a real $3\times 3$ matrix $S$ that describes correlations between Alice and Bob. 

The effect of these correlations is manifest when Bob performs a quantum measurement on his qubit. Bob's measurement is described by a set of $n$ operators $\{E_i\}$ satisfying $E_i\geq 0$ and $\sum_{i=1}^n E_i=\mathbb{1}$~\cite{NC}. We write $E_i=M_i^\dagger M_i$; on obtaining measurement outcome $E_i$, the two-qubit state $\rho$ is \textit{collapsed} according to $\rho\rightarrow(\mathbb{1}\otimes M_i) \rho (\mathbb{1}\otimes M_i^\dagger)/\tr(\rho (\mathbb{1}\otimes E_i))$. When this measurement-induced collapse occurs, the Bloch vector representing Alice's qubit is changed. By choosing which measurement to perform, Bob can thus \textit{steer} Alice. Given all possible measurements by Bob, the set of Bloch vectors to which Alice can be steered forms an ellipsoid inside the Bloch sphere  $\widetilde\B$~\cite{QSE}. This is Alice's steering ellipsoid, which we denote $\E$. 

The centre of $\E$ is given by $\vec d$. Since $\widetilde\B$ is centred at the origin, the distance between the centres of $\E$ and $\widetilde\B$ is $d=|\vec d|$. The eigenvalues of the real $3 \times 3$ symmetric matrix $Q=SS^\mathrm{T}$ correspond to the squared lengths of the ellipsoid semiaxes; the eigenvectors of $Q$ give the orientation of these axes. That $\E$ is contained inside $\widetilde\B$ is necessary but not sufficient to meet the constraint $\rho\geq 0$. Note also that $\E$ could be a degenerate ellipsoid (i.e. an ellipse, line or point), corresponding to $Q$ being rank deficient.

\vspace{-3mm}
\section{Formulation of an inequality for the existence of a nested tetrahedron}

Steering ellipsoids give a new geometric perspective on two-qubit entanglement in the form of the \textit{nested tetrahedron condition}~\cite{QSE}. Any two-qubit state $\rho$ that is not entangled is described as separable and may be decomposed into $m\leq 4$ terms as $\rho=\sum_{i=1}^m p_i \alpha_i \otimes \beta_i$, where $\alpha_i, \beta_i$ are single qubit states and $\{p_i\}$ is a probability distribution~\cite{Wootters}. This leads to the result that $\rho$ is separable if and only if there exists a tetrahedron circumscribed about $\E$ and inscribed in $\widetilde\B$. For the case that $\E$ is an ellipse, Poncelet's porism can be used to show that if there exists a tetrahedron nested between $\E$ and $\widetilde\B$ then there must also exist a triangle nested between $\E$ and $\widetilde\B$~\cite{Poncelet}. 

Independently of the nested tetrahedron condition we can also formulate an algebraic condition for separability. Define $\rho^\Gamma=(\mathbb{1} \otimes \mathrm{T})(\rho)$, where $\mathrm{T}$ is the usual matrix transposition map. For a two-qubit state we write $\rho=\big(\begin{smallmatrix} A&B \\ C&D\end{smallmatrix}\big)$, where $A,B,C$ and $D$ are $2\times 2$ matrices; then $\rho^\Gamma=\left(\begin{smallmatrix} A^\mathrm{T}&B^\mathrm{T} \\ C^\mathrm{T}&D^\mathrm{T}\end{smallmatrix}\right)$. Using the Peres-Horodecki criterion~\cite{Peres,Horodecki}, Ref.~\cite{QSE} shows that the condition $\det \rho^\Gamma\geq 0$ is necessary and sufficient for separability\footnote{For example, take $\rho_1=\frac{1}{4}\Bigg(\begin{smallmatrix} 1&0&0&0 \\ 0&1&1&0 \\ 0&1&1&0 \\ 0&0&0&1 \end{smallmatrix}\Bigg)$ and $\rho_2=\frac{1}{2}\Bigg(\begin{smallmatrix} 0&0&0&0 \\ 0&1&1&0 \\ 0&1&1&0 \\ 0&0&0&0 \end{smallmatrix}\Bigg)$ -- as required, these are unit trace and positive semi-definite. We find $\rho_1^\Gamma=\frac{1}{4}\Bigg(\begin{smallmatrix} 1&0&0&1 \\ 0&1&0&0 \\ 0&0&1&0 \\ 1&0&0&1 \end{smallmatrix}\Bigg)$ and $\rho_2^\Gamma=\frac{1}{2}\Bigg(\begin{smallmatrix} 0&0&0&1 \\ 0&1&0&0 \\ 0&0&1&0 \\ 1&0&0&0 \end{smallmatrix}\Bigg)$. Computing the signs of $\det \rho_1^\Gamma$ and $\det \rho_2^\Gamma$ then immediately identifies $\rho_1$ as separable and $\rho_2$ as entangled.}. For a canonical state of the form \eqref{canonical}, this determinant can easily be evaluated in terms of $Q$ and $\vec d$~\cite{QSE}. A canonical state $\rho$ is seen to be separable if and only if
\begin{equation}\label{sep}d^4-2 \widetilde u d^2 +\widetilde q\geq 0,\end{equation}
where $\widetilde u =1-\tr Q+2\vec{\hat d}^\mathrm{T}Q\vec{\hat d}$ and $\widetilde q =1-2\tr Q-8 \sqrt{\det Q} +2\tr(Q^2)- (\tr Q)^2$, with unit vector $\vec{\hat d} = \vec d / d$.

We thus have two equivalent necessary and sufficient conditions for the separability of $\rho$: the nested tetrahedron condition and inequality \eqref{sep}. It is not at all obvious how either one of these conditions could be found directly from the other. Crucially, however, we now have an algebraic formulation for the existence of a nested tetrahedron. Converting the unit sphere $\widetilde\B$ to a sphere $\B$ of radius $R$, we arrive at the key inequality that will be used to derive all our results in 2 and 3 dimensions.
\begin{thm}\label{keyinequality}Let $\E$ be an ellipsoid, described by matrix $Q$, contained inside a sphere $\B$ of radius $R$. The eigenvalues of $Q$ give the squared lengths of the semiaxes of $\E$; the eigenvectors give the orientation of the axes. The centre of $\E$ relative to the centre of $\B$ is $\vec d$. The sole necessary and sufficient condition for the existence of a tetrahedron circumscribed about $\E$ and inscribed in $\B$ is
\begin{equation}\nonumber d^4-2 u d^2 +q\geq 0,\end{equation}
where $u =R^2-\tr Q+2\vec{\hat d}^\mathrm{T}Q\vec{\hat d}$ and $q =R^4-2R^2\tr Q -8 R \sqrt{\det Q}+2\tr(Q^2)- (\tr Q)^2$, with unit vector $\vec{\hat d} = \vec d / d$. For a degenerate $\E$, the tetrahedron can always be taken to be a triangle.
\begin{proof}
Inequality \eqref{sep} is necessary and sufficient for the existence of a tetrahedron circumscribed about $\E$ and inscribed in $\widetilde\B$. The result then immediately follows from scaling the quantities $\widetilde u$ and $\widetilde q$ to $u$ and $q$ respectively.
\end{proof}
\end{thm}

Note that only the term involving $\vec{\hat d}^\mathrm{T}Q\vec{\hat d}$ depends on the orientation of $\vec d$ rather than just its magnitude $d$. This \emph{skew} term describes how the axes of an ellipsoid are oriented relative to the centre vector $\vec d$.

\vspace{-3mm}
\section{Derivation of the Grace-Danielsson inequality}

Theorem \ref{keyinequality} is all that will be needed to derive inequality \eqref{2d_ineq} for the 2-dimensional scenario and the corresponding 3-dimensional Grace-Danielsson inequality \eqref{3d_ineq}. We will also look at an example of how Theorem \ref{keyinequality} can be used to find new results for ellipses and ellipsoids.

\begin{cor}\label{circle}Let $\E$ be a circle of radius $r$ contained inside a circle $\B$ of radius $R$. The distance between the centres of $\E$ and $\B$ is $d$. The sole necessary and sufficient condition for the existence of a triangle circumscribed about $\E$ and inscribed in $\B$ is \begin{equation}\nonumber d^2\leq R(R-2r),\end{equation}
giving the result \eqref{2d_ineq}.
\begin{proof}$\E$ may be described by $Q = \textrm{diag}(r^2, r^2, 0)$. The degenerate case of Theorem \ref{keyinequality} gives a condition for the existence of a triangle circumscribed about $\E$ and inscribed in a sphere of radius $R$. Setting $\vec d=(d_1, d_2, 0)$ ensures that $\E$ and $\vec d$ are coplanar, so that Theorem \ref{keyinequality} equivalently gives a condition for the existence of a triangle circumscribed about $\E$ and inscribed in a circle of radius $R$. The skew term $\vec{\hat d}^\mathrm{T}Q\vec{\hat d}=r^2$ does not depend on the orientation of $\vec d$ within its plane. Evaluating $u$ and $q$ gives $u=R^2$ and $q=R^4-4R^2r^2$. Theorem \ref{keyinequality} then tells us that a nested triangle exists if and only if $d^4-2 u d^2 +q=(d^2-R^2)^2-4R^2r^2 \geq 0$, from which the result follows.
\end{proof}
\end{cor}

The physical significance of this result in terms of steering ellipsoids inside the Bloch sphere is considered in detail in Ref.~\cite{MaximalQSE}. Briefly, the `no pancake theorem' states that the equatorial plane of the Bloch sphere is not a valid steering ellipsoid as it does not achieve $\rho \geq 0$; by setting $R=1$, Corollary \ref{circle} extends this result to identify the largest circular $\E$ that is a valid steering ellipsoid for a given $\vec d$.

\begin{cor}\label{sphere}Let $\E$ be a sphere of radius $r$ contained inside a sphere $\B$ of radius $R$. The distance between the centres of $\E$ and $\B$ is $d$. The sole necessary and sufficient condition for the existence of a tetrahedron circumscribed about $\E$ and inscribed in $\B$ is \begin{equation}\nonumber d^2\leq (R+r)(R-3r),\end{equation}
giving the Grace-Danielsson result \eqref{3d_ineq}.
\begin{proof}$\E$ may be described by $Q = \textrm{diag}(r^2, r^2, r^2)$. The skew term $\vec{\hat d}^\mathrm{T}Q\vec{\hat d}=r^2$ does not depend on the orientation of $\vec d$. Evaluating $u$ and $q$ gives $u=R^2-r^2$ and $q=R^4-6R^2r^2-8Rr^3-3r^4$. Theorem \ref{keyinequality} then tells us that a nested tetrahedron exists if and only if $d^4-2 u d^2 +q=(d-R-r)(d+R+r)(d^2-(R+r)(R-3r))\geq 0$, from which the result follows.
\end{proof}
\end{cor}

Again, there is a physical interpretation of this result in quantum information theory~\cite{MaximalQSE}. By identifying the largest spherical $\E$ that corresponds to a separable state $\rho$, the entangled-separable state boundary can be found for a special class of $\rho$ known as inept states~\cite{Inept}.

Finally, we consider results for an ellipsoid $\E$. Note that Theorem \ref{keyinequality} concerns \emph{any} ellipse or ellipsoid inside a sphere -- in particular, this includes ellipsoids that are oriented with no semiaxis collinear with $\vec d$. Such ellipsoids have an awkward skew term $\vec{\hat d}^\mathrm{T}Q\vec{\hat d}$, which causes significant difficulties when attempting to formulate algebraic conditions for when $\E$ is contained inside a sphere~\cite{EinB}. Remarkably, Theorem \ref{keyinequality} works in full generality to give conditions for the existence of a nested tetrahedron for an ellipsoid with any skew. The degenerate case can be used to give conditions for the existence of a nested triangle for any ellipse inside a circle or sphere. We believe this to be the first formulation of necessary and sufficient conditions for these general scenarios.

As an example of how Theorem \ref{keyinequality} can be used for ellipsoids, we give a result for a specially oriented class of ellipsoid.

\begin{cor}\label{ellipsoid}Let $\E$ be an ellipsoid with semiaxes $s_1,s_2,s_3$ contained inside a sphere $\B$ of radius $R$. The distance between the centres of $\E$ and $\B$ is $d$, and $\E$ has its $s_1$ axis collinear with $\vec d$. The sole necessary and sufficient condition for the existence of a tetrahedron circumscribed about $\E$ and inscribed in $\B$ is \begin{equation}\nonumber d^2\leq (R-s_1)^2-(s_2+s_3)^2.\end{equation}
\begin{proof}$\E$ may be described by $Q = \textrm{diag}(s_1^2, s_2^2, s_3^2)$ and $\vec d=(d, 0, 0)$. Evaluating $u$ and $q$ gives $u=R^2+s_1^2-s_2^2-s_3^2$ and $q=R^4-2R^2(s_1^2+s_2^2+s_3^2)-8Rs_1s_2s_3+s_1^4+s_2^4+s_3^4-2s_1^2s_2^2-2s_2^2s_3^2-2s_3^2s_1^2$. Theorem \ref{keyinequality} then tells us that a nested tetrahedron exists if and only if $d^4-2 u d^2 +q=(d^2-R^2-s_1^2+s_2^2+s_3^2)^2-(2Rs_1+2s_2 s_3)^2\geq 0$, from which the result follows.
\end{proof}
\end{cor}

Note that by setting $s_1=s_2=s_3=r$ Corollary \ref{ellipsoid} reproduces the result for spheres given in Corollary \ref{sphere}.

\vspace{-3mm}
\section{Conclusions}

Using the quantum steering ellipsoid formalism and an understanding of the separability of two-qubit states, we have given a necessary and sufficient condition for the existence of a tetrahedron circumscribed about an ellipsoid and inscribed in a sphere. This condition can be used to immediately derive a number of results in 2- and 3-dimensional Euclidean geometry, most notably the Grace-Danielsson inequality for a nested tetrahedron. As well as providing an entirely new formulation for the general case of an ellipsoid, our results are of interest for the novelty of their derivations.

\subsection*{Acknowledgements}

The author acknowledges useful discussions with Mowaffaq Hajja, David Jennings, Sania Jevtic, Matthew Pusey and Terry Rudolph. This work was supported by EPSRC.

\end{document}